\documentclass[11pt]{article}

\usepackage{amsmath,amssymb,amsthm}
\usepackage{geometry}
\newtheorem{theorem}{Theorem}[section]
\newtheorem{lemma}[theorem]{Lemma}
\newtheorem{proposition}[theorem]{Proposition}

\theoremstyle{definition}
\newtheorem{definition}[theorem]{Definition}
\newtheorem{remark}[theorem]{Remark}

\newtheorem{conj}[theorem]{Conjecture}
\newtheorem{que}[theorem]{Question}

\title{Finite Groups with Few  Almost Maximal Subgroups\footnote{Corresponding author: W. Meng (mlwhappyhappy@163.com).
W. Meng is supported by National Natural Science Foundation of China (12661001), Guangxi Natural Science
				Foundation Program (2026GXNSFAA00641090) and Guangxi Colleges and Universities Key Laboratory of Data Analysis and Computation.
}}
\author{Wei Meng\textsuperscript{a},
Boru Zhang\textsuperscript{b}, Jiakuan Lu\textsuperscript{b}\\
\textsuperscript{a}School of Mathematics and Computing Science, Guilin University of Electronic Technology,\\
 Guilin, Guangxi, 541002, P.R. China.\\
  \textsuperscript{b}School of Mathematics and Statistics, Guangxi Normal University,\\
   Guilin, Guangxi, 541004, P.R. China.\\
   Emails: mlwhappyhappy@163.com, Brzhangyq@163.com, jklu@gxnu.edu.cn}

\date{}

\begin{document}

\maketitle

\begin{abstract}
Let $G$ be a finite group. L. Morgana, C. E. Praeger,   A.E. Zalesski \cite {MOR} introduced the notion of almost maximal subgroup of $G$. A subgroup $H$ of $G$ is said to be  almost maximal in $G$ if $N_G(H)$ is maximal in $G$. In this paper, let $AM(G)$ denote the number of conjugacy classes of all almost maximal subgroups of $G$. First, we prove that $G$ is Dedekindian if and only if $AM(G)=0$. Second, we show that the finite groups with $AM(G)\leq 1$ are supersolvable. Furthermore, we characterize the structure of finite groups with $AM(G)=1$. Final, we investigate a lower bound of $AM(G)$ for certain non-solvable groups.
\end{abstract}



\section{Introduction}

Throughout this paper, all groups are finite.
Since the completion of the Classification of Finite Simple Groups, one of the central problems in finite group theory has been the determination of the maximal subgroups of almost simple groups, namely, groups $X$ such that $S\unlhd X\leq
Aut(S)$ for some finite nonabelian simple group $S$ (see, for instance, the discussion and references in \cite [ Sections 3 and 4] {KLE} and the recent survey \cite{RON}). For alternating and symmetric groups $X \in \{A_n,S_n\}$, M. W. Liebeck, C. E. Praeger and J. Saxl \cite {LIE} reduced the problem  to the classification of all inclusions $H<L$  in which both $H$ and $L$  are almost simple primitive permutation groups of degree $n$ and $N_L(H)$ is maximal in $L$. In addition, chains
\begin{center}
$H_1 <H_2 <\cdots<H_{r+1} <X$
\end{center}
 of arbitrary length $r$ were constructed in \cite[pp. 368-369]{LIE} so that each inclusion $H_i<H_{i+1}$ has these properties.

More recently,  L. Morgana, C. E. Praeger,   A.E. Zalesski \cite {MOR}  continued the investigation of maximal subgroups of finite alternating and symmetric groups. Rather than fixing  $X\in\{A_n,S_n\}$, however, they turned their attention to an arbitrary finite nonabelian simple group $T$. Generally,  $T$ admits many embeddings $T<X$ with
$X \in\{A_n,S_n\}$ for which $T$ acts primitively on the underlying set$\Omega=\{1,2,\ldots,n\}$, such embeddings are called primitive embeddings of
 $T$ in $X$.  This action is equivalent to the right multiplication action of $T$
on the coset space $\{Mx| x\in  T\}$ of a maximal subgroup $M<T$, where $n = |T : M|$.

 L. Morgana, C. E. Praeger and A. E. Zalesski \cite {MOR} directed their attention to those embeddings for which the normaliser
$N_X(T)$  is a maximal subgroup of $X$, designating such embeddings as '\textbf{almost maximal}'.

\begin{definition} \cite {MOR}
A subgroup $T$ of a group $X$ is said to be \textbf{almost maximal} in $X$ if $N_X(T)$
is a maximal proper subgroup of $X$.
\end{definition}

Let  $T$ be a finite simple group and $X$ be a finite alternating group
or a finite symmetric group. L. Morgana, C. E. Praeger and A. E. Zalesski \cite {MOR} established that, if
$T$ is almost maximal in $X$ and $N_X(T)$  is not primitive, then $T$  is necessarily either cyclic of order at most
$3$ or an alternating group. Conversely, they proved that, apart from the single exception $PSL(2,11)$,
 every finite nonabelian simple group $T$  occurs for some $n$  as an almost maximal subgroup of an alternating group
$A_n$ in such a way that $T$  acts primitively; furthermore, the primitive embedding can be chosen so that
$N_{A_n}(T)$   is the unique maximal subgroup of $A_n$ containing $T$.\\

The above works have focused primarily on almost simple groups, especially on maximal subgroups of alternating and symmetric groups and on almost maximal embeddings, with emphasis on determining when $N_T(X)$ is maximal and on the primitivity of the corresponding action.
For solvable groups, however, almost maximal subgroups have not yet been systematically studied, whether with respect to their structure, the distribution of their conjugacy classes, or their influence on the group. This gap motivates us to turn from almost simple groups to solvable groups and to investigate how the number  of conjugacy classes of almost maximal subgroups reflects and restricts the structure of a solvable group. \\

  For convenience, let $G$ be a finite group and let
  $AM(G)$ denote the number of conjugacy classes of almost maximal subgroups of $G$.  In this paper, we investigate the structure of finite groups under the assumption that $AM(G)$  is small.\\

  Recall that a group $G$ is called a Dedekind group if   every subgroup of $G$ is normal in $G$.  First, we provide a new characterization of Dedekind groups in terms of $AM(G)$.

\begin{theorem}
Let $G$ be finite group. Then $AM(G)=0$ if and only if  $G$ is a Dedekind group.
\end{theorem}

As a natural next step, we consider the case where $AM(G)\leq 1$  and obtain the following supersolvability criterion.

\begin{theorem}
Let $G$ be a finite group. If $AM(G)\leq 1$, then $G$ is supersolvable.
\end{theorem}

Furthermore, we character the finite group with $AM(G)=1$.

\begin{theorem}
Let $G$ be a finite group with $AM(G)=1$. Then one of the following statements holds.

(1) If $G$ is nilpotent,  then $G$ is a non-Dedekind $p$-group.

(2) If $G$ is non-nilpotent, then  $|\pi(G)|=2$.
\end{theorem}

\begin{theorem}
Suppose that $G$ is a finite non-nilpotent group  with  $AM(G)=1$. Then   $G=P\rtimes Q$,   $P\in Syl_p(G)$ is cyclic, and $Q\in Syl_q(G)$ is Dedekindian and $|Q/Q_G|=q$, where $p$ and $q$ are distinct primes, and $q\mid p-1$. In particular, if $|P|=p$, then $Q$ is also cyclic.
\end{theorem}

All unexplained notations and terminologies
    are standard and can be found in \cite {DO, HB}. In addition, $\pi(G)$ denotes the set of primes dividing
the order of $G$; $C_m$ denotes the cyclic group of order $m$; $C_{p}^d$  denotes the elementary abelian $p$-group
of order $p^d$; $A \times B$ means the direct product of $A $ and $B$. $ A \rtimes B$ means a split extension of the normal
subgroup $A$ by the non-normal subgroup $B$; $M_G=\bigcap\limits_{g\in G}M^g$ denotes the core of $M$ in $G$; $[M]$ denotes the $G$-conjugacy class of $M$.
\section{Preliminaries}

We begin with some elementary observations.

\begin{lemma}
Let $G$ be a finite group and $N$ a normal subgroup of $G$. Then $AM(G/N)\leq AM(G)$.
\end{lemma}

{\noindent\bf  Proof.} First, observe that $N_{G/N}(H/N)=N_G(H)/N$ for any subgroup  $H/N$ of $G/N$. It follows that $H/N$ is almost maximal in $G/N$ if and only if $H$ is almost maximal in $G$.   Furthermore, $H/N$ and $K/N$ are conjugate in $G/N$ if and only if $H$ and $K$ are conjugate in $G$. Therefore the correspondence $H/N$  induces an injection from the conjugacy classes of almost maximal subgroups of $G/N$
 into those of $G$   and so $AM(G/N)\leq AM(G)$. $\Box$

\begin{lemma}
Let $G$ be a finite group and $M$ a maximal subgroup of $G$. Then $M$ is almost maximal in $G$ if and only if $M=N_G(M)$, i.e., $M$ is non-normal in $G$.
\end{lemma}

{\noindent\bf  Proof.} It follows directly from the definition of   almost maximal subgroup. $\Box$

\begin{lemma} \cite {BEL2}
Suppose that a $p$-group $G$ is not Dedekindian. Then there is $H<G$ such that $|G:N_G(H)|=p$ and $|H|\geq |G'|$.
\end{lemma}

\begin{lemma} \cite {BEL1}
If a finite group $G$ has at most $2$ conjugacy classes of non-normal maximal subgroups,
then $G$ is solvable.
\end{lemma}

\begin{lemma} \cite {BEL1}
Let $G$ be a finite non-solvable group. Then $G$ has three conjuagacy classes of maximal
subgroups if and only if either $G/\Phi(G)\cong PSL(2, 7)$ or $PSL(2, 2^p)$, where $p$ is a prime.
\end{lemma}

\begin{lemma} \cite {HB}
Let $G$ be a finite group. Then $G$ is supersolvable if and only if $|G:M|$ is a prime for any  maximal subgroup $M$ of $G$.
\end{lemma}

\begin{lemma}
Let $G$ be a finite group. If $AM(G)\leq 2$, then $G$ is solvable.
\end{lemma}

{\noindent\bf  Proof.}  Let $t$ be the number of conjugacy classes of non-normal maximal subgroups of $G$. Then $G$ has at least $t$ conjugacy classes of almost maximal subgroups by Lemma 2.2. By hypothesis, $AM(G)\leq 2$ implies that $t\leq 2$. Applying Lemma 2.4, $G$ is solvable. $\Box$

\begin{lemma}
Let $G=P\rtimes Q$, and $P\cong C_p^d$$(d\ge 2)$,  $Q\cong C_q$,  
where $p$ and $q$ are distinct primes, and $q\mid p-1$. Suppose that $G$ has exactly one conjugacy class
of non-normal maximal subgroups. Then the unique non-normal maximal subgroup
is isomorphic to ${C_p^{\,d-1}\times C_q}$.  
\end{lemma}

 {\noindent\bf  Proof.}
Since $q\mid p-1$, then $Q$ acts on $P$ by conjugation which is complete reduction by Maschke theorem.  So we can assume that 
\begin{center}
$P=V_1\times V_2\times\cdots\times V_d$,
\end{center}
where each $V_i\cong C_p$ is $Q$-invariant.   
Put $W_i=\prod\limits_{j\neq i}V_j$, then $W_i$ is $Q$-invariant maximal subgroup of $P$.  Thus, $M_i=W_iQ$ is maximal subgroup of $G$ of index $p$.

We claim that $M_i=W_iQ\lhd G$ if and only if $[P,Q]\leq W_i$. 
Indeed,  $W_iQ\lhd G$ is equivalent to $Q^x\leq W_iQ$ for every $x\in P$. 
Since $u^x=[x,u^{-1}]u$, we have $ Q^x=[x,Q]Q$. Therefore, $[x,Q]Q\subseteq W_iQ$ if and only   $[x,Q]\subseteq W$ for every $x\in P$.
Hence,  $[P,Q]\le W_i$. 

Moreover, observe that $P=[P,Q]\times C_P(Q)$ which implies that  $M_i=W_iQ$ is normal in $G$ if and only if $V_i\leq C_P(Q)$, i.e., $Q$ acts trivially on $V_i$.  In other words, $M_i$ is non-normal in $G$ if and only if $Q$ acts on-trivially on $V_i$.

Let $i\ne j$, and suppose that $Q$ acts non-trivially on both $V_i$ and
$V_j$. Assume that $M_i^g=M_j,\qquad g\in G$. Then 
\begin{center}
$W_i=W_{i}^g=(M_i\cap P)^g=M_{i}^g\cap P=M_j\cap P=W_j$
\end{center} 
 as $W_i$ is normal in $G$.  
It follows that the number of conjugacy classes of non-normal maximal
subgroups of $G$ is exactly $r=log_{p}^{|[P,Q]|}$.  By hypothesis, $G$ exactly one conjugacy class of non-normal maximal
subgroups. Hence $|[P,Q]|=p$. Without loss of generality, assume that $[P,Q]=V_1$, then $C_P(Q)=W_1$. Consequently, $M_1=W_1Q\cong C_{p}^{d-1}\times C_q$ is a representation of the unique conjugacy class of non-normal maximal subgroups 
 of $G$. $\Box$

\section{The case $AM(G)=0$}

In this section, we   prove Theorem 1.2.\\

 {\noindent\bf  Proof of Theorem 1.2.} First, suppose that $G$ is a Dedekind group. Then  $H$   is  normal in $G$, i.e., $N_G(H)=G$ for every proper subgroup $H$ of $G$. This implies that $AM(G)=0$.

 Conversely, suppose that $AM(G)=0$. For any maximal subgroup $M$ of $G$, $M$ is not almost maximal in $G$. It follows that $M$ is normal in $G$ by Lemma 2.2. So $G$ is nilpotent. If $G$ is not a Dedekind group, there exists a Sylow-p-subgroup $P$ of $G$ which is non-Dedekindian for some $p\in\pi(G)$. Appling Lemma 2.3, $P$ has a subgroup $H$ such that $|P:N_P(H)|=p$. Since $G$ is nilpotent,  we get $|G:N_G(H)|=|P:N_P(H)|=p$. Thus $N_G(H)$ is a maximal subgroup of $G$. It follows that $H$ is almost maximal in $G$, a contradiction. So $G$ is Dedekindian. The proof of theorem is complete. $\Box$\\

 \section {The case  $AM(G)=1$}

In this section, we focus on finite groups satisfying  $AM(G)=1$. We first state the following proposition.

 \begin{proposition}
Let $G$ be a non-nilpotent    group with $AM(G)=1$ and $M$ be   a maximal subgroup of $G$. Then every non-trivial normal subgroup of $M$ is normal in $G$.
\end{proposition}

{\noindent\bf  Proof.}  Suppose that $G$ is a non-nilpotent group. Then $G$ has at least  one non-normal maximal subgroup, say, $H$.  By Lemma 2.2, $H$ is almost maximal in $G$.  $AM(G)=1$ implies that $H$ is a representation of  the unique conjugacy classes of almost maximal subgroups of $G$. For any maximal subgroup $M$ of $G$, let $N$ be a proper normal subgroup of $M$. Then $M\leq N_G(N)\leq G$. If $N_G(N)=M$, then $N$ is almost maximal in $G$. However, $N\neq H^g$ for any $g\in G$. This leads to $AM(G)\geq2$, a contradiction. Therefore, $N_G(N)=G$ and hence $N$ is normal in $G$.
  $\Box$ \\

 Next, we will give the proof of Theorem 1.3.\\

  {\noindent\bf  Proof of Theorem 1.3.}  By  Theorem 1.2,  we may assume that $AM(G)=1$. Furthermore,  if all maximal subgroup of $G$ are normal in $G$, then $G$ is nilpotent and hence $G$ is also  supersolvable.   So we   assume that $G$ has a non-normal maximal subgroup, say, $M$. By Lemma 2.2, $M$ is almost maximal in $G$.  Moreover,   $AM(G)=1$   implies that the set  $\{M^g\mid g\in G\}$ is   the   unique conjugacy classes of  almost maximal subgroups of $G$, and $\{M^g\mid g\in G\}$ is also the unique conjugacy class of non-normal maximal subgroups of $G$.  It follows that for any maximal subgroup $K$ of $G$ with $K\not\in\{M^g\mid g\in G\}$, $K$ is normal in $G$.  By Lemma 2.7, we know that  $G$ is solvable, and thus $|G:K|$ is a prime number for   any maximal subgroup $K$ of $G$ with $K\not\in\{M^g\mid g\in G\}$.\\

   By Lemma 2.6, to prove that $G$  is supersolvable, we need only show that $|G:M|$ is also a prime number. We now proceed by induction on
$|G|$.  Let $M_G$ denote the core of $M$ in $G$. Our proof will be divided into the following two cases.\\

Case I: $M_G\neq 1$.

 In this case, consider quotient group $G/M_G$. Obvious, $M/M_G$ is  a non-normal maximal subgroup of  $G/M_G$, and hence $M/M_G$ is almost maximal in $G/M_G$ by Lemma 2.2.   It follows that $AM(G/M_G)\neq 0$. On the other hand, by Lemma 2.1, we have $AM(G/M_G)\leq AM(G)=1$ which obliges to $AM(G/M_G)=1$. So $G/M_G$ satisfies the theorem hypothesis. By induction, we get  that $G/M_G$ is supersolvable.  Since $M/M_G$ is a maximal subgroup of $G/M_G$, $|G/M_G:M/M_G|$ is a prime number by Lemma 2.6. Consequently, $|G:M|=|G/M_G:M/M_G|$ is a prime number.  \\

 Case II: $M_G=1$.

 Let $M'$ be  the commutator subgroup of  $M$.  Observe that $G$ is solvable, hence $M$ is also solvable. Consequently,  $M'$ is a proper normal subgroup of $M$. Applying Proposition  4.1,  $M'$ is also normal in $G$ which implies $M'\leq M_G=1$.  Thus,  $M$ is abelian.

  We claim that $M$ is group of order prime. Otherwise,  $M$ has at least one  non-trivial subgroup, say $L$.  Then $L$ is normal in $M$ as $M$ is abelian. Applying Proposition 4.1 again, $L$ is normal in $G$. Thus $L\leq H_G=1$. This is absurd.  The claim as desired.

In the following,  let $N$ be a minimal normal subgroup of $G$.  The solvability of $G$ implies that  $N$ is an elementary abelian group.  Moreover, Since $M$ is a group of order prime, we have $N\not\leq M$ and $N\cap M=1$. By the maximality of $M$, we get that  $G=NM$. It follows that $|G:N|=|M|$ is a prime number and hence $N$ is a maximal subgroup of $G$.

Moreover,  we show that $N$ is also a group of order prime. If not, then $N$ has a non-trivial subgroup, say $H$.  Since $N$ is an elementary abelian group, $H$  is   normal in $N$. By Proposition 4.1, $H$ is normal in $G$. Therefore, $HM$ is a subgroup of $G$. By maximality of $M$, we have $G=HM$. It follows that $H=N$, a contradiction. So $N$ is a group of order prime.    Consequently, $|G:M|=|N|$ is a prime number.
  The proof of theorem is complete. $\Box$\\

We now proceed to prove Theorem 1.4.\\

  {\noindent\bf  Proof of Theorem 1.4.} Let $G$ be a finite group with $AM(G)=1$.
  First, assume that $G$ is nilpotent. By Theorem 1.2, $G$ is not Dedekindian.
Therefore, there exists a   Sylow $p$-subgroup $P$ of $G$  for some $p\in\pi(G)$ such that $P$ is not Dedekindian. Applying Lemma 2.3, $P$ has a proper subgroup $L$ such that $|P:N_P(L)|=p$. Suppose that $G$ has a non-trivial  normal Hall $p'$-subgroup $H$. Then  $G=P\times H$. Moreover, we get  $N_G(L)=N_G(LH)=N_P(L)\times H$ is a maximal subgroup of $G$. So both $L$ and $LH$ are almost maximal in $G$. This deduces that $AM(G)\geq2$, a contradiction. So $G$ is a $p$-group. The conclusion (1) holds.\\

In what follows, assume that $G$ is non-nilpotent. Let $M$ be a non-normal maximal subgroup of $G$. Then $[M]$ is  the unique conjugacy class of almost maximal subgroups of $G$.  We   proceed by induction on
$|G|$.  Since $AM(G)=1$, the group $G$ is supersolvable by Theorem 1.3.
Let $N$ be a minimal normal subgroup of $G$. Then $N$ has   prime
order, say $|N|=r$.  By Lemma 2.1, we have $AM(G/N)\leq AM(G)=1$.  Our proof will be divided into the following two cases.\\

Case I: $G/N$ is nilpotent.

We claim that $N\nleq M$. Otherwise, $M/N$ would be a
maximal subgroup of the nilpotent group $G/N$, and
therefore $M/N$ is  normal in $G/N$. This imply
$M\trianglelefteq G$, a contradiction. So $N\nleq M$.  Consequently,
  $G=NM$ and $N\cap M=1$.

 If $AM(G/N)=1$, then $G/N$ is a $p$-group for some prime $p$ by conclusion (1).  Since $G$ is non-nilpotent, we have $p\neq r$ and hence $|\pi(G)|=2$. The conclusion (2) holds.

 Now, we assume that $AM(G/N)=0$. Applying Theorem 1.2, $M\cong G/N$ is a Dedkind group.
 By Proposition 4.1, every non-trivial subgroup $L$ of $M$ is normal in $G$. Thus $NL=N\times L$
Observe that $M$ does not centralize $N$, there exists  an element $x\in M$ such that $x\not\in C_G(N)$. So   $M=\langle x\rangle$ is a cyclic group. Furthermore, if $|\pi(M)|\geq 2$, then    Then the cyclic group $M$ has two distinct
maximal subgroups $M_1,M_2$ satisfying $M=M_1M_2$.  By above arguments, we know that
both $M_1$ and $M_2$ centralize $N$, so $M$ centralizes
$N$, a contradiction. Therefore, $M$ is a cylic $p$-group for some $p\neq r$. It follows that $\pi(G)=\{p,r\}$. The conclusion (2) holds again.\\

 Case II: $G/N$ is not nilpotent.

In this case,   $AM(G/N)=1$ and hence $G/N$ satisfies the theorem hypothesis. By induction, we get that $|\pi(G/N)|=2$.  Suppose that $\pi(G/N)=\{p,q\}$.
Suppose  that $r\notin\pi(G/N)=\{p,q\}$.  Then $N$ is a normal Sylow $r$-subgroup of $G$.
By the Schur-Zassenhaus theorem, $N$ has a complement
$H$ in $G$, so
\begin{center}
$G=N\rtimes H$, and  $H\cong G/N$.
\end{center}
Observe that  $|G:H|=|N|=r$ is a prime, then $H$ is a maximal subgroup of $G$.

We claim that $N\leq Z(G)$. If not, then $H$ is a non-normal maximal subgroup of $G$ and thus $H$ is almost maximal in $G$ by Lemma 2.2.
On the other hand, since $G/N$ is non-nilpotent, $G/N$ has a non-normal maximal subgroup, say $M/N$. It follows that $M$ is a non-normal maximal subgroup of $G$  and hence $M$ is almost maximal in $G$. However, $M$ cannot be conjugate to $H$ as $N\leq M$ and $N\cap H=1$. This implies that $AM(G)\geq 2$, a contradiction. So the claim as desired.

Now, we have $N\leq Z(G)$. Consequently, $G=N\times H$ and $AM(H)=AM(G/N)=1$.   There exists a proper subgroup $L$ of $H$
 such that $N_H(L)$ is maximal in $H$. Moreover,
 \begin{center}
 $N_G(L)=N_G(NL)=N\times N_H(L)$
  \end{center}
  is maximal in $G$.  This shows that both $L$ and $NL$ are almost maximal in $G$ which contradicts to $AM(G)=1$. Therefore, $r\in\pi(G/N)$ and so $|\pi(G)|=2$.  The conclusion (2) is obtained.
 The proof of theorem is complete.
$\Box$\\

To further investigate the structure of finite groups $G$  with $AM(G)=1$,   we prove the following proposition.

   \begin{proposition}
Assume that $G$ is a  non-nilpotent group with   $AM(G)=1$,
and  $M$ is a non-normal maximal subgroup of $G$.  Then
for every normal subgroup $N\trianglelefteq G$, we have
\begin{align*}
AM(G/N)=
\begin{cases}
1,&\text{if }N\leq M_G,\\[2mm]
0,&\text{if }N\nleq M_G.
\end{cases}
\end{align*}
\end{proposition}

{\noindent\bf  Proof.}
By Lemma 2.1, $AM(G/N)\leq AM(G)=1$.
If $N\leq M_G$, then $N\leq M$. Since $M$ is non-normal in $G$, $M/N$ is non-normal in $G/N$. So $M/N$ is almost in $G/N$. It follows that $AM(G/N)=1$.

 Suppose that $N\nleq M_G$ and $AM(G/N)=1$. Let $H/N$ be an almost maximal subgroup of $G/N$. Consequently, $H$ is almost maximal in $G$. So there is an element $g\in G$ such that $H=M^g$. This leads to $N\leq M_{G}^g=M_G$, a contradiction. Therefore, $AM(G/N)=0$ if $N\nleq M_G$.
 $\Box$\\

 \begin{proposition}
 Let $G$ be a non-nilpotent group  and $M$ be a non-normal maximal subgroup of $G$. If $AM(G)=1$, then $G/M_G\cong C_q\rtimes C_p$, where $p$ and $q$ are distinct primes, and $p\mid q-1$.
 \end{proposition}

{\noindent\bf  Proof.} Suppose that $G$ is a non-nilpotent group with $AM(G)=1$. By Theorem 1.4 (2),   we may assume that $\pi(G)=\{p,q\}$, where $p<q$. Let $M$ be a non-normal maximal subgroup of $G$ and $M_G$ the core of $M$ in $G$. Consider the quotient group $\overline{G}=G/M_G$. Then $\overline{M}=M/M_G$ is a non-normal maximal subgroup of $\overline{G}$ and hence $\overline{M}$ is core-free in $\overline{G}$. Moreover, by Proposition 4.2, we have $AM(\overline{G})=1$. By Proposition 4.1, we know that every non-trivial normal subgroup of $\overline{M}$ is normal in $\overline{G}$. It follows that the commutator subgroup $\overline{M}'$ of $\overline{M}$ is normal in $\overline{G}$. So $\overline{M}'\leq \overline{M}_{\overline{G}}=1$, and hence $\overline{M}$ is abelian. Applying Proposition 5.1 again, $\overline{M}\cong C_p$ is a cyclic group of order $p$.

On the other hand, by Theorem 1.3, we know that $G$ is supersolvable and so is $\overline{G}$. Let $\overline{N}$ be a minimal normal subgroup of $\overline{G}$. Then $\overline{N}\not\leq\overline{M}$ and hence $\overline{N}\cong C_q$. Consequently, $\overline{G}=\overline{N}\rtimes\overline{M}\cong C_q\rtimes C_p$ and $p\mid q-1$. $\Box$\\

{\noindent\bf  Proof of Theorem 1.5.} Suppose  that $G$ is a non-nilpotent group with $AM(G)=1$. Then $G$ possesses a unique conjugacy class of non-normal maximal subgroups, say $[M]$.  Then $|G:M|=p$.   By Theorem 1.4(2), we can let $\pi(G)=\{p,q\}$$(p>q)$. Let $M$ be a non-normal maximal subgroup of $G$. Applying Proposition 4.3,  $G/M_G\cong C_p\rtimes C_q$ and  $q\mid p-1$.

Let $P\in Syl_p(G)$ and $Q\in Syl_q(G)$. By Theorem 1.3,  $G$ is supersolvable  and hence  $G=P\rtimes Q$.   Proposition By
 Since $PM_G/M_G$ is a Sylow $p$-subgroup of $G/M_G$, we get  $P\nleq M_G$. Applying Proposition 4.2, we have $AM(Q)=AM(G/P)=0$ and so $Q$ is a Dedekind group. Moreover, $M/M_G =QM_G/M_G\cong Q/Q\cap M_G$ which implies that $|Q/Q_G|=q$. 
 
 In what follows, we will prove that $P$ is cyclic.  Suppose that $P$ is non-cyclic. Then $P/\Phi(P)\cong C_{p}^d$, where $d\geq2$. Let $N=\Phi(P)Q_G$. Then $N\leq M_G$.  Consider the  quotient group $\overline{G}=G/N$, then $\overline{G}$ is also non-nilpotent.  Furthermore,  by Proposition 4.2. $AM(G)=1$. It follows that $[\overline{M}]$ is the   unique conjugacy class of non-normal maximal subgroups of $\overline{G}$. Hence, $\overline{M}$ is almost maximal in $\overline{G}$. 
 
  On the other hand, 
 \begin{center}
 $\overline{G}=G/N=PQ/N=PN/N\rtimes QN/N\cong P/\Phi(P)\rtimes Q/Q_G\cong C_{p}^d\rtimes C_q$.
 \end{center}
 By Lemma 2.8, $\overline{M}\cong C_{p}^{d-1}\times C_q$. This implies that $\overline{Q}$ is normal in $\overline {M}$ and hence $\overline{Q}$ is almost maximal in $\overline{G}$. Since $AM(\overline{G})=1$, we get that $\overline{Q}=\overline{M}$. So $d=1$, a contradiction.  
 
 Final, suppose that $|P|=p$, then $Q=M$ is a representation of the unique conjugacy class of non-normal maximal subgroups of $G$. By above arguments, $Q$ is a Dedekind group. Applying Proposition 4.1, every proper subgroup of $Q$ is normal in $G$. Therefore, if $Q$ is non-cyclic, then $Q$ will be normal in $G$. This is a contradiction. 
 The proof of theorem is complete. $\Box$

 \begin{remark}
 Let $G$ be a finite group with $AM(G)=2$. Then $G$ may be not supersolvable. For example,
 let $G=A_4$ be an alternating group of degree $4$.  It is easy to check that the subgroups of $G$ of order
$2$   and $3$ are all almost maximal in  $G$; Therefore, $AM(G)=2$. However, $G$ is not a supersolvable group.
 \end{remark}

\section{Some finite groups $G$  with small $AM(G)$}
In this section, we first give some examples of finite groups $G$ with $AM(G)=1$. 

\begin{theorem}
Let $G=D_{2n}$ be a dihedral group of order $2n$. Then $AM(G)=1$ if and only if $n=p^{\alpha}$, where $p$ is an odd prime.
\end{theorem}
{\noindent\bf  Proof.} Let $G=\langle a,b\mid a^n=b^2=1,b^{-1}ab=a^{-1}\rangle$ be a dihedral group of order $2n$. If $n=2$, then $G$ is abelian and hence $AM(G)=0$. So we may assume that $n\geq3$.

First, assume that $AM(G)=1$. By Theorem 1.5, $|\pi(G)|\leq 2$. Suppose that $n=2^{\alpha}$, where $\alpha\geq2$. Write
$M_1=\langle a^2,b\rangle$ and $M_2=\langle a^2,ab\rangle$. Then each $M_i$ is non-cyclic maximal subgroup of $G$ for $i=1,2$. Let $H_1=\langle a^4,b\rangle$ and $H_2=\langle a^4,ab\rangle$. Then $H_i\lhd M_i$. But $H_i$ is non-normal in $G$. So $N_G(H_i)=M_i$. Since $M_i$ is normal in $G$, $H_1$ cannot be conjugate to $H_2$. It follows that $AM(G)\geq2$, a contradiction. Thus, $n=p^{\alpha}$ is an odd prime power.

Conversely, assume that $n=p^{\alpha}$ is an odd prime power. Then maximal subgroups of $G$ forms  two conjugacy classes as follows, $[\langle a\rangle]$ and $[\langle a^p,b\rangle]$. Moreover, every subgroup of $\langle a\rangle$ is normal in $G$ and $\langle a^p,b\rangle$ is  a non-normal maximal subgroup of $G$. Thus, $\langle a^p,b\rangle$ is almost maximal in $G$ by Lemma 2.2. It is easy to see that $AM(G)=1$ if $\alpha=1$. So we can let $\alpha\geq2$.   Then every non-trivial normal subgroup $\langle a^{p^m}\rangle$ of  $\langle a^p,b\rangle$ is normal in $G$,  where $1\leq m\leq \alpha-1$. This implies that $AM(G)=1$. The proof of theorem is complete.
$\Box$

\begin{theorem}
Let $G=Q_{4n}$ be a generalized quaternion group of order $4n$, where $n\geq2$. Then $AM(G)=1$ if and only if $n=p^{\alpha}$, where $p$ is an odd prime.
\end{theorem}
 {\noindent\bf  Proof.} Suppose that $G=\langle a,b\mid a^{2n}=1,a^n=b^2,b^{-1}ab=a^{-1}\rangle$. Then $Z(G)=\langle a^n\rangle\cong C_2$ and $G/Z(G)\cong D_{2n}$.  Applying Theorem 5.1, we can get the conclusion of theorem. The proof of theorem is complete. $\Box$\\

By Theorem 1.4(1), we know that if $G$ is a nilpotent group with $AM(G)=1$, then $G$ is a $p$-group. Brandl \cite {BRA}  characterized the finite $p$-groups   whose non-normal subgroups are all conjugate. Such $p$-groups are as follows 
\begin{center}
$M_p=\langle a,b\mid a^{p^{n-1}}=b^p,a^b=a^{1+p^{n-2}}\rangle$, 
\end{center}
where $p$ is a prime and $n\geq 3$ if $p \geq3$, or $n \geq4$ if $p =2$. Applying Lemma 2.3, we have $AM(M_P)=1$. The following question would seem to merit some consideration.

\begin{que}
Classified the finite $p$-groups with $AM(G)=1$. 
\end{que}

Next, we investigate a lower bound for the number of conjugacy classes of almost maximal subgroups in non-solvable groups. 

 \begin{proposition}
 (1) $AM(A_5)=6$;

 (2) $AM(PSL(2,7))=8$;

 (3) $AM(PSL(2,2^p))\geq 6$.
 \end{proposition}

 {\noindent\bf  Proof.} (1) Let $G=A_5$ be an alternating group of degree $5$. The maximal subgroups of $G$ fall into three conjugacy classes:
$[A_4]$, $[D_{10}]$, and $[S_3]$. Since  $G=A_5$ is a simple group, every maximal subgroup $M$ of $G$ is  self-normalizing.
  Hence these three maximal subgroups are almost maximal and contribute three conjugacy classes.

Certain proper subgroups of these maximal subgroups also have maximal normalizers. Namely, $N_G(V_4)=A_4$, $N_G(C_5)=D_{10}$,
 $N_G(C_3)=S_3$, where $V_4\leq A_4$, $C_5\leq D_{10}$ and  $C_3\leq S_3$.
 Hence, each of these subgroups is almost maximal and contributes one conjugacy class.
Therefore  $AM(G)=6$.

(2) Let $G = PSL(2, 7)$. Then $G$ has exactly three conjugacy classes of maximal subgroups, say, $[S_4]_1$, $[S_4]_2$, and $[C_7\rtimes C_3]$.
Hence these three maximal subgroups are almost maximal and contribute three conjugacy classes. Furthermore, observe that $C_7$ is normal in $C_7\rtimes C_3$,  $A_4$ and $V_4$ are normal in $S_4$. Thus each of these subgroups is almost maximal and contributes one conjugacy class. Therefore, $AM(G)=8$.

(3) Let $G=PSL(2,q)$, where $q=2^p$, $p$ is prime.   By Dickson's classification \cite [Hauptsatz II.8.27] {HB}, for
 the maximal subgroups of $G$
 fall into exactly three conjugacy classes:

(a) Borel subgroups $B\cong C_{2}^p\rtimes C_{q-1}$;

(b) Dihedral subgroups $D_{2(q+1)}$;

(c) Dihedral subgroups $D_{2(q-1)}$.

Every member of above maximal subgroups    has at least one  non-trivial normal subgroup,  and thus they  are almost maximal in $G$. This implies that $AM(G)\geq6$. $\Box$ 

\begin{theorem}
Let $G$ be a finite non-solvable group. If $G$ has three conjuagacy classes of maximal
subgroups, then $AM(G)\geq 6$. 
\end{theorem}

 {\noindent\bf  Proof.}  By Lemma 2.5, $G/\Phi(G)\cong PSL(2,7)$ or $PSL(2,2^p)$, where $p$ is a prime. Applying Proposition 5.3, $AM(G/\Phi(G))\geq6$. So we get that $AM(G)\geq AM(G/\Phi(G))\geq 6$ by Lemma 2.1. $\Box$\\

 Recall that a group $G$ is called a minimal non-solvable group if $G$ is non-solvable, but every proper subgroup of $G$ is solvable.  By result, if $G$ is a minimal non-solvable group, then $G/\Phi(G)$ is a minimal simple group, i.e., $G/\Phi(G)$ is a simple group all of whose proper subgroups are solvable. Thompson \cite {THO} determined all minimal simple groups.

\begin{theorem}
Let $G$ be a minimal non-solvable group. Then $AM(G)\geq6$.
\end{theorem}
 {\noindent\bf  Proof.} Let $G$ be a minimal non-solvable group. Consider the quotient group $\overline{G}=G/\Phi(G)$, then $\overline{G}$ is a minimal simple group. Let $\overline{M_1}$, $\overline{M_2}$,$\ldots$, $\overline{M_t}$ be  representations of   conjugacy classes of all maximal subgroups of $\overline{G}$. Then $t\geq3$ by Lemma 2.4. Applying Lemma 2.2, every $\overline{M_i}$ is almost maximal in $\overline{G}$. Furthermore, since each $\overline{M_i}$ is solvable, $\overline{M_i}$ has at least one non-trivial normal subgroup, say $\overline{H_i}$. Consequently,  $N_{\overline{G}}(\overline{H_i})=\overline{M_i}$ and hence $\overline{H_i}$ is almost maximal in $\overline{G}$. This yields to $AM(\overline{G})\geq 2t\geq 6$.  Finial, applying Lemma 2.1, we have $AM(G)\geq AM(\overline{G})\geq 6$. The proof of theorem is complete.  $\Box$\\

 We end this paper by proposing the following conjecture.

 \begin{conj}
Every finite group $G$ satisfying $AM(G)\leq 5$ is solvable.
\end{conj}

\end{document}